\documentclass[10pt]{article}
\usepackage{amsfonts}
\usepackage{amsmath}
\usepackage{authblk}
\usepackage{hyperref}

\normalbaselines

\title{Controlling a nonlinear Fokker-Planck equation via inputs with nonlocal action}
\author{\c{S}tefana-Lucia Ani\c{t}a \\
\small ``Octav Mayer'' Institute of Mathematics of the Romanian Academy, Ia\c{s}i 700506, E-mail: stefi$\underline{~}$anita@yahoo.com} 
\date{}

\begin{document}

\begin{titlepage}
\maketitle
\thispagestyle{empty}

\abstract{This paper concerns an optimal control problem $(P)$ related to a nonlinear Fokker-Planck equation. The problem is deeply related to a stochastic 
optimal control problem $(P_S)$ for a McKean-Vlasov equation. 
The existence of an optimal control is obtained for the deterministic problem $(P)$. 
The existence of an optimal control is established and necessary optimality conditions are derived for a penalized optimal control problem $(P_h)$ related to a backward Euler approximation of the nonlinear Fokker-Planck equation (with a constant discretization step $h$).
Passing to the limit ($h\rightarrow 0$) one derives the necessary optimality conditions for problem $(P)$.
}
\vspace{3mm}

\vspace{2mm}


\end{titlepage}

\section{Formulation of the problem}

\large

Consider the following optimal control problem:
$$\underset{\zeta \in {\cal M}}{\text{Minimize}} \  \int_0^T \int_{\mathbb{R}^d}G(t,x)\rho ^{\zeta}(t,x) dx \ dt+\int_{\mathbb{R}^d}G_T(x)\rho ^{\zeta}(T,x)dx+\int_{\mathbb{R}^d}Q(x,\zeta (x))dx,\leqno({\bf P})$$
where $T\in (0,+\infty )$, $d\in \mathbb{N}, \ d\geq 3$, and $\rho ^{\zeta}$ is the weak solution to the following nonlinear Fokker-Planck equation:

\begin{equation}\label{nfp1}
\left\{ \begin{array}{ll}
\displaystyle {{\partial \rho}\over {\partial t}}(t,x)=-\nabla \cdot (u(x)b(\rho (t,x))\rho (t,x))+\Delta \beta (\rho (t,x)), & t\in [0,T] , \ x\in \mathbb{R}^d \\
~ & ~ \\
\rho (0,x)=\rho _0(x), & x\in \mathbb{R}^d .
                                     \end{array}
                          \right. 
\end{equation}
Here $u(x)=-\nabla \left( \int_{\mathbb{R}^d}{\cal G}_R(x-x')\zeta (x')dx'\right)={\cal K}(\zeta )(x)$, 
$$\zeta \in {\cal M}=\{ \theta \in L^{\infty }(\mathbb{R}^d); \ 0\leq \theta (x)\leq \tilde{M}_0 \ \mbox{\rm a.e. } x\in \mathbb{R}^d, \ \theta (x)=0 \ \mbox{\rm a.e. } |x|_d>R_0\};$$
$\tilde{M}_0$, $R$, $R_0$ are positive constants, and we assume that ${\cal G}_R\in C_0^2(\mathbb{R}^d)$, ${\cal G}_R(x)>0$ if $|x|_d<R$ and ${\cal G}_R(x)=0$ if $|x|_d\geq R$.
We have that ${\cal K}:L^{\infty }(\mathbb{R}^d)\longrightarrow L^{\infty }(\mathbb{R}^d)^d$ is defined by 
$${\cal K}(\zeta )(x)=-\nabla ({\cal G}_R*\zeta )(x)$$
(where ${\cal G}_R*\zeta $ is the convolution product of ${\cal G}_R$ and $\zeta $).

The weak solution to (\ref{nfp1}) may be viewed as the probabilistic density of a population while $\zeta (x)$ is the density at $x \in \mathbb{R}^d$ of a second population (or of another entity, which may be a substance or a signal) which produces a stimulus to the first population. 
This second population is time-independent, imobile, located in $\overline{B(0_d;R_0)}$ $(R_0>0$), and it repels the individuals of the first population which are at a distance less than $R$. It means that $\zeta $ is an input (control) with nonlocal action. This action, denoted by $u$, is expressed mathematically in terms of the so-caled ``generalized gradient'' (nonlocal gradient) with kernel ${\cal G}_R$.

Actually, $u(x)=-\nabla ({\cal G}_R*\zeta )(x)={\cal K}(\zeta )(x)$ 
describes the nonlocal action (effect) of $\zeta$ (the second population) towards the individuals of the first population located at $x\in \mathbb{R}^d$. 
The term $-\nabla \cdot (u(x)b(\rho (t,x))\rho (t,x))=\nabla \cdot (\nabla ({\cal G}_R*\zeta )(x)b(\rho (t,x))\rho (t,x))$ describes a cross-diffusion (see e.g. \cite{capasso, KM}), and leads to a particular form of Keller-Segel's models (see e.g. \cite{burger1, burger2}).

Problem $(P)$ asks to optimally displace a population via the repeller action produced by a second population of density $\zeta $. It is a natural requirement (due to logistic constraints) to place this second   
population in a bounded domain (here we have considered it to be a ball) and to have bounded density.
\vspace{2mm}

Notice that the case of a local repeller action at $x\in \mathbb{R}^d$ is a limit case of the nonlocal one, where
$${\cal G}_R=\delta $$
($\delta $ is the Dirac distribution) and
$$u(x)=-\nabla \zeta (x).$$

If the second population attracts the first one, then we have to take $u(x)=\nabla ({\cal G}_R*\zeta )(x)$ instead (and $u(x)=\nabla \zeta (x)$ in the case of local action). 
\vspace{5mm}

Assume that the following hypotheses hold:
\begin{itemize}
\item[{\bf (H1)}] $\beta \in C^2(\mathbb{R})$, $\beta (0)=0$ and $\exists 0<\gamma _0<\gamma _1$ such that
$$\gamma _0\leq \beta '(r)\leq \gamma _1, \ \forall r\in \mathbb{R},\quad  \mbox{\rm and } \Psi \in C^2(\mathbb{R}), \ \mbox{\rm where }\Psi (r)={{\beta (r)}\over r};$$
\item[{\bf (H2)}] $b\in C^1(\mathbb{R})$, $b$ is bounded, $b(r)\geq 0$, $\forall r\geq 0$ and 
$$|b^*(r)-b^*(\tilde{r})| \leq \alpha | \beta (r)-\beta (\tilde{r})|, \quad \forall r, \tilde{r} \in \mathbb{R},$$ where $\alpha >0$ and $b^*(r)=b(r)r$;
\item[{\bf (H3)}] $\rho _0\in   L^1(\mathbb{R}^d)\cap  L^{\infty }(\mathbb{R}^d)$, $\rho _0(x)\geq 0$ a.e. $x\in \mathbb{R}^d$, $\int_{\mathbb{R}^d}\rho _0(x)dx=1$;
\item[{\bf (H4)}] $G \in C_b([0,T]\times \mathbb{R}^d)\cap L^2((0,T)\times \mathbb{R}^d)$, $G_T\in C_b(\mathbb{R}^d) \cap H^1(\mathbb{R}^d)$, 
$Q:\mathbb{R}^d\times \mathbb{R}\longrightarrow \mathbb{R}$, $Q|_{B(0_d;R_0)\times [0,\tilde{M}_0]}\in C_b^{0,1}(B(0_d;R_0)\times [0,\tilde{M}_0])$,
$$0\leq G(t,x), \ \forall (t,x)\in [0,T]\times \mathbb{R}^d, \quad 0\leq G_T(x), \ \forall x\in \mathbb{R}^d,$$
$$0\leq Q(x,z), \ \forall (x,z)\in \mathbb{R}^d\times [0,\tilde{M}_0], \quad Q(x,z)=0, \ \forall x\in \mathbb{R}^d, |x|_d>R_0, \forall z\in [0,\tilde{M}_0],$$
and the mapings $z \mapsto Q(x,z)$ are convex with respect to $z\in [0,\tilde{M}_0]$ for any $x\in \mathbb{R}^d$.
\end{itemize}

The following notations $\nabla =\left({{\partial }\over {\partial x_1}}, {{\partial }\over {\partial x_2}}, \cdots , {{\partial }\over {\partial x_d}}\right)$, 
$\Delta = \sum _{i=1}^d {{\partial ^2}\over {\partial x_i ^2}}$ will be used. Sometimes we will use the simplified notation 
$\| \cdot \| _{L^k}$ (instead of $\| \cdot \| _{L^k(\mathbb{R}^d)}$ or $\| \cdot \| _{L^k(\mathbb{R}^d;\mathbb{R}^d)}$) for the usual norm of $L^k(\mathbb{R}^d)$ or of $L^k(\mathbb{R}^d;\mathbb{R}^d)$.
\vspace{2mm}

For some recent results concerning the optimal control of Fokker-Planck equations we refer to \cite{anita3, AB1, AB2, barbu3, FG}.
The difficulty of our present study is mainly due to the nonlinearities in the Fokker-Planck equation (\ref{nfp1}). 
\vspace{3mm}

In section 2 we will recall that there exists a unique weak solution to the nonlinear Fokker-Planck equation (\ref{nfp1}) (see \cite{barbu3}) and that our optimal control problem is deeply related to a stochastic optimal control problem $(P_S)$ related to a McKean-Vlasov equation. 
In section 3 one proves the existence of an optimal control for problem $(P)$. The existence of an optimal control is established and necessary optimality conditions are derived for a penalized optimal control problem $(P_h)$ related to a backward Euler scheme for the nonlinear Fokker-Planck equation (with a constant discretization step $h$) in section 4. Using the optimality conditions for $(P_h)$ one derives the necessary optimality conditions for problem $(P)$ in section 5.
Section 6 contains final remarks and comments.

\section{The nonlinear Fokker-Planck equation and its relationship to a McKean-Vlasov equation}

For any $\zeta \in {\cal M}$ we consider the equation (\ref{nfp1}) as the following Cauchy problem in $L^1(\mathbb{R}^d)$:
\begin{equation}\label{eq2.8}
\left\{ \begin{array}{ll}
\displaystyle {{d\rho }\over {dt}}(t)+A^{\zeta}\rho (t)=0, \quad t\in [0,T] \\
~ \\
\rho (0)=\rho _0,
                                     \end{array}
                          \right. 
\end{equation}
where $A^{\zeta}:D(A^{\zeta})\subset L^1(\mathbb{R}^d)\longrightarrow L^1(\mathbb{R}^d)$ is an m-accretive operator in $L^1(\mathbb{R}^d)$ given by 
$$D(A^{\zeta })=\{ \theta \in L^1(\mathbb{R}^d); \ -\Delta \beta (\theta )+\nabla \cdot ({\cal K}(\zeta )b(\theta )\theta ) \in L^1(\mathbb{R}^d)\} ,$$
$$A^{\zeta}(\theta ) = -\Delta \beta (\theta )+\nabla \cdot ({\cal K}(\zeta )b(\theta )\theta ) \quad \mbox{\rm in } {\cal D}' (\mathbb{R}^d ), \quad \theta  \in D(A^{\zeta })$$
(see Proposition 2.3 in \cite{BR5}). Moreover, $\overline{D(A^{\zeta})}=L^1(\mathbb{R}^d)$ and $(I+\lambda A^{\zeta})^{-1}({\cal P})\subset {\cal P}, \ \forall \lambda >0$, where ${\cal P}$ is the set of all probability densities on $\mathbb{R}^d$.
\vspace{2mm}

\noindent
\bf Definition 2.1. \rm
Let $\zeta \in {\cal M}$. A function $\rho ^{\zeta} \in C([0,T];L^1(\mathbb{R}^d))$ is called a weak/mild solution to (\ref{nfp1})/(\ref{eq2.8}) if
$$\rho ^{\zeta}(t)=\underset{h\rightarrow 0+}{\text{lim}} \tilde{\rho }^{\zeta}_h(t) \ \mbox{\rm in } L^1(\mathbb{R}^d), \ \mbox{\rm uniformly for } t \in [0,T],$$
where for any $h>0$, $\tilde{\rho }^{\zeta}_h:[0,T]\rightarrow L^1(\mathbb{R}^d)$ is given by
$$\tilde{\rho }^{\zeta }_h(t)=\rho _h^{\zeta ,i} \quad \mbox{\rm if } t\in [ih, (i+1)h), t\leq T, \ i=0,1,...,N=\left[{T\over h}\right]$$
and $\rho ^{\zeta ,i+1}_h$ is the solution to 
\begin{equation}\label{eq2.7}
\left\{ \begin{array}{ll}
\rho ^{\zeta ,0}_h=\rho _0, \quad \rho ^{\zeta ,i}_h\in L^1(\mathbb{R}^d), &  i=1,...,N, \\
~ & ~ \\
\rho ^{\zeta ,i+1}_h-h\Delta \beta (\rho ^{\zeta ,i+1}_h)+h\nabla \cdot ( {\cal K}(\zeta ) b(\rho ^{\zeta ,i+1}_h)\rho ^{\zeta ,i+1}_h)=\rho ^{\zeta ,i}_h \ \mbox{\rm in } {\cal D}'(\mathbb{R}^d), & i=0,1,...,N-1.
                                     \end{array}
                          \right. 
\end{equation}
\vspace{2mm}

The existence and uniqueness of such a weak/mild solution $\rho ^{\zeta }$ to (\ref{nfp1})/(\ref{eq2.8}) follow via Theorem 2.1 in \cite{BR5} (which is based on Crandall-Liggett's existence theorem (see \cite{barbu1})). 
Moreover, by Theorem 2.1 in \cite{BR5} we get that $\rho ^{\zeta }$ is also a distributional solution to (\ref{nfp1}). Furthermore, we have that $\rho ^{\zeta}(t) \in {\cal P}, \ \forall t \in [0,T]$ and $\rho ^{\zeta} \in L^{\infty}((0,T)\times \mathbb{R}^d)$.
\vspace{2mm}

\noindent
\bf Remark 2.1. \rm If we consider $N\in \mathbb{N}^*$ and $h={T\over N}$, then we get that
$$\rho ^{\zeta}(t)=\underset{h\rightarrow 0+}{\text{lim}} \rho ^{\zeta}_h(t)=\underset{N\rightarrow +\infty }{\text{lim}} \rho ^{\zeta}_h(t) \ \mbox{\rm in } L^1(\mathbb{R}^d), \ \mbox{\rm uniformly for } t \in [0,T],$$
where $\rho ^{\zeta}_h:[0,T]\rightarrow L^1(\mathbb{R}^d)$ is given by
$$\rho ^{\zeta}_h(t)=\left\{ \begin{array}{ll}
\rho _0, \quad & \mbox{\rm if } t=0, \\
~ & ~ \\
\rho ^{\zeta,i+1}_h, \quad & \mbox{\rm if } t\in (ih, (i+1)h], \quad i=0,1,...,N-1.
                                                   \end{array}
                                         \right. $$

Indeed, if $t=ih$ ($i\in \{ 0,1,..., N\}$), then
$\rho _h^{\zeta }(t)=\rho _h^{\zeta ,i}=\tilde{\rho }_h^{\zeta }(t)$
and so
$$\| \rho _h^{\zeta }(t)-\rho ^{\zeta }(t)\| _{L^1(\mathbb{R}^d)}=\| \tilde{\rho }_h^{\zeta }(t)-\rho ^{\zeta }(t)\| _{L^1(\mathbb{R}^d)}.$$
On the other hand, if $t\in (ih,(i+1)h), \ i\in \{ 0,1,...,N-1\} ,$ then
$\rho _h^{\zeta }(t)=\rho _h^{\zeta }((i+1)h)$ and so
$$\| \rho _h^{\zeta }(t)-\rho ^{\zeta }(t)\| _{L^1(\mathbb{R}^d)}=\| \rho _h^{\zeta }((i+1)h)-\rho ^{\zeta }(t)\| _{L^1(\mathbb{R}^d)}=\| \tilde{\rho }_h^{\zeta }((i+1)h)-\rho ^{\zeta }(t)\| _{L^1(\mathbb{R}^d)}$$
$$\leq \| \tilde{\rho }_h^{\zeta }((i+1)h)-\rho ^{\zeta }((i+1)h)\| _{L^1(\mathbb{R}^d)}+\| \rho ^{\zeta }((i+1)h)-\rho ^{\zeta }(t)\| _{L^1(\mathbb{R}^d)}.$$
Taking into account the convergence of $\tilde{\rho }^{\zeta }_h(t)$ to $\rho ^{\zeta }(t)$ in $L^1(\mathbb{R}^d)$, uniformly for $t\in [0,T]$ and the fact that
$\rho ^{\zeta }\in C([0,T];L^1(\mathbb{R}^d))$, we get the conclusion. 
\vspace{5mm}

Consider now the following controlled McKean-Vlasov SDE on $\mathbb{R}^d$:
 
\begin{equation}\label{eq1.1}
\left\{ \begin{array}{ll}
dX(t)=f\left(X(t),{\cal K}(\zeta )(X(t)), {{d{\cal L}_{X(t)}}\over{dx}}(X(t))\right)dt+\sigma \left(X(t),{{d{\cal L}_{X(t)}}\over{dx}}(X(t))\right)dW(t), \quad t\in[0,T] \\
~ \\
X(0)=X_0,
                                     \end{array}
                          \right. 
\end{equation}

\noindent
where ${\cal L}_Y$ is the law/distribution of the random variable $Y$ and ${{d{\cal L}_Y}\over{dx}}$ is its density (when it exists) with respect to the Lebesgue measure $dx$,
and
$$f(x,u,r)=ub(r), \quad \sigma (x,r)=\left( {{2\beta (r)}\over r}\right)^{1\over 2}.$$
Since for any $\zeta \in {\cal M}$, $\rho ^{\zeta }$ is a distributional solution to the nonlinear Fokker-Planck solution, which is also $t$-narrowly continuous, then by superposition's principle 
(see \cite{figalli, trevisan, HRW}) we conclude that there exists a unique (in law) probabilistically weak solution $(X^{\zeta }(t))_{t\in [0,T]}$ to (\ref{eq1.1}) (which has $\rho ^{\zeta }(t)$ as a probability density).
Actually, the weak solution $X^{\zeta }$ corresponds to a probability space $(\Omega, {\cal F}, \mathbb{P})$ with normal filtration $({\cal F}_t)_{t \in [0,T]}$ 
and $(W(t))_{t \in [0,T]}$ is an $\mathbb{R}^d-({\cal F}_t)_{t \in [0,T]}$ Brownian motion.

Hence, 
$$ \mathbb{E}\left[ \int_0^TG(t,X^{\zeta }(t))dt\right] +\mathbb{E}[G_T(X^{\zeta }(T))]= \int_0^T \int_{\mathbb{R}^d}G(t,x)\rho ^{\zeta }(t,x) dx \ dt+\int_{\mathbb{R}^d}G_T(x)\rho ^{\zeta }(T,x)dx.$$
It follows that the next stochastic optimal control problem
$$\underset{\zeta \in {\cal M}}{\text{Minimize}} \  \mathbb{E}\left[ \int_0^TG(t,X^{\zeta }(t))dt\right]+\mathbb{E}[G_T(X^{\zeta }(T))]+\int_{\mathbb{R}^d}Q(x,\zeta (x))dx,\leqno({\bf P_S})$$
is equivalent to the deterministic optimal control problem $(P)$.
\vspace{3mm}

We mention that lately there is a high interest in McKean-Vlasov SDEs and in stochastic optimal control problems. Some recent and important results concerning the McKean-Vlasov 
can be found in \cite{barbu3,BR2,BR3,BR4}. As concerns problem $(P_S)$ (which is a stochastic optimal control problem related to McKean-Vlasov SDE) we see that it is deeply related to
an optimal control problem $(P)$ for a nonlinear Fokker-Planck equation. A similar approach, reducing the investigation of a stochastic optimal control problem to a deterministic optimal control problem, 
has been recently used for the case of feedback controllers independent of time in \cite{anita,anita2,barbu2} and in \cite{anita3,FG} for controllers dependent of time as well, for optimal control problems with other cost functionals and related to stochastic differential equations with drift and diffusion independent of the probabilistic density. 
For other approaches for stochastic optimal control problems see \cite{BBD,BBT,BRZ} while for standard results concerning the control of stochastic differential equations we refer to the monographs \cite{FR,oksendal}.

\section{Existence of an optimal control for problem $\mbox{\bf (P)}$}

Let $m^*=\inf _{\zeta \in {\cal M}} I(\zeta )\in [0,+\infty )$, where
$$I(\zeta )=\int_0^T\int_{\mathbb{R}^d}G(t,x)\rho ^{\zeta }(t,x)dx \ dt+\int_{\mathbb{R}^d}G_T(x)\rho ^{\zeta }(T,x)dx+\int_{\mathbb{R}^d}Q(x,\zeta (x))dx.$$
For any $\zeta \in {\cal M}$ we have that ${\cal K}(\zeta )\in C_0^1(\mathbb{R}^d;\mathbb{R}^d)$,
$$\begin{array}{ll}
supp \ {\cal K}(\zeta )\subset \overline{B(0_d;R+R_0)} , \\
\displaystyle
|{\cal K}(\zeta )(x)|_d\leq \tilde{M}_0\int_{\mathbb{R}^d}|\nabla {\cal G}_R(-x')|_d dx'=M, \quad \forall x\in \mathbb{R}^d , \\
|\nabla {\cal K}(\zeta )(x)|_{d\times d}\leq M_1, \quad \forall x\in \mathbb{R}^d,
\end{array}$$
where $M_1$ is a constant independent of $\zeta \in {\cal M}$.

For any sequence $\{ \zeta _k\} _{k\in \mathbb{N}^*}\subset {\cal M}$, there exists a subsequence $\{ \zeta _{k_l}\}_{l\in \mathbb{N}^*}$ and $\zeta ^*\in {\cal M}$
such that
$$\zeta _{k_l}\longrightarrow \zeta ^* \quad \mbox{\rm weakly * in } L^{\infty }(\mathbb{R}^d).$$
This implies that
$${\cal K}(\zeta _{k_l})(x)=-\int_{\mathbb{R}^d}\nabla {\cal G}_R(x-x')\zeta _{k_l}(x')dx' \longrightarrow -\int_{\mathbb{R}^d}\nabla {\cal G}_R(x-x')\zeta ^*(x')dx'={\cal K}(\zeta ^*)(x), \quad \forall x\in \mathbb{R}^d,$$
and that
$$\nabla {\cal K}(\zeta _{k_l})(x)\longrightarrow \nabla {\cal K}(\zeta ^*)(x), \quad \forall x\in \mathbb{R}^d.$$
Let us notice that by Arzel\`{a}'s compactness theorem we may infer that there exists a subsequence $\{ {\cal K}(\zeta _{k_{l_r}})\}_{r\in \mathbb{N}^*}$ such that
$${\cal K}({\zeta _{k_{l_r}}})\longrightarrow \tilde{u}, \quad \mbox{\rm uniformly for } x\in \overline{B(0_d;R+R_0)}.$$
Since ${\cal K}( \zeta_{k_{l_r}})(x)=0, \ \forall |x|_d>R+R_0$, we get that
$${\cal K}({\zeta _{k_{l_r}}})\longrightarrow \tilde{u }, \quad \mbox{\rm uniformly for } x\in  \mathbb{R}^d$$
(where $\tilde{u}$ has been extended by the value $0_d$ outside $ \overline{B(0_d;R+R_0)}$). Since ${\cal K}({\zeta _{k_{l_r}}})\longrightarrow {\cal K}({\zeta ^*})$,
$\forall x\in \mathbb{R}^d$, we conclude that
$${\cal K}({\zeta _{k_{l_r}}})\longrightarrow {\cal K}({\zeta ^*}), \quad \mbox{\rm uniformly for } x\in  \mathbb{R}^d$$
\vspace{3mm}

We postpone for the time being proof of the next auxiliary result:
\vspace{3mm}

\noindent
\bf Lemma 3.1. \it If $\{ \zeta_n \}_{n\in \mathbb{N}^*} \subset {\cal M}$, $\zeta _n \longrightarrow \zeta \in {\cal M}$ weakly * in $L^{\infty}(\mathbb{R}^d)$,
$${\cal K}(\zeta _n)=-\nabla ({\cal G}_R*\zeta _n)\longrightarrow -\nabla ({\cal G}_R*\zeta )={\cal K}(\zeta ), \quad \mbox{\it uniformly for } x\in \mathbb{R}^d,$$
then for any $h>0$ sufficiently small and for any $f\in L^1(\mathbb{R}^d)$ we have that
$$y_n \longrightarrow y \quad \mbox{\it in } L^1(\mathbb{R}^d),$$
where $y_n, \ y\in L^1(\mathbb{R}^d)$,
\begin{equation}\label{eq3.2}
y_n-h\Delta \beta (y_n)+h\nabla \cdot ( {\cal K}(\zeta _n) b^*(y_n))=f \quad \mbox{\it in } {\cal D}'(\mathbb{R}^d)
\end{equation}
and
\begin{equation}\label{eq3.3}
y-h\Delta \beta (y)+h\nabla \cdot ( {\cal K}(\zeta ) b^*(y))=f \quad \mbox{\it in } {\cal D}'(\mathbb{R}^d).
\end{equation}
\rm
\vspace{2mm}

\noindent
\bf Remark 3.1. \rm By Lemma 2.5 in \cite{BR5} we have that equations (\ref{eq3.2}) and (\ref{eq3.3}) have unique solutions, respectively.
Moreover, if one uses more accurate notations: $y_n(f)$ for the solution to (\ref{eq3.2}), and $y(f)$ for the solution to  (\ref{eq3.3}), then
$$\| y_n(f)-y_n(g)\|_{L^1(\mathbb{R}^d)}\leq \|f-g\| _{L^1(\mathbb{R}^d)},$$
for any $f, \ g\in L^1(\mathbb{R}^d)$. By Lemma 3.1 one obtains that
$$\| y(f)-y(g)\|_{L^1(\mathbb{R}^d)}\leq \|f-g\| _{L^1(\mathbb{R}^d)}, \quad f, \ g\in L^1(\mathbb{R}^d).$$

\vspace{3mm}

Using now Trotter-Kato's theorem (Theorem 2.1-p.241, \cite{barbu}) we may conclude that
\vspace{3mm}

\noindent
\bf Lemma 3.2. \it If $\{ \zeta_n \}_{n\in \mathbb{N}^*} \subset {\cal M}$, $\zeta _n \longrightarrow \zeta \in {\cal M}$ weakly * in $L^{\infty}(\mathbb{R}^d)$,
${\cal K}(\zeta _n)\longrightarrow {\cal K}(\zeta ),$  uniformly for $x\in \mathbb{R}^d,$
then
$$\rho^{\zeta _n} \longrightarrow \rho^{\zeta } \quad \mbox{\it in } C([0,T];L^1(\mathbb{R}^d)).$$
\vspace{3mm}

\noindent
\bf Theorem 3.3. \it Problem $(P)$ has at least one optimal control. 
\vspace{2mm}

\noindent
Proof. \rm Let
$\{ \zeta _k\} _{k\in \mathbb{N}^*}\subset {\cal M}$ such that
$$m^*\leq I(\zeta _k)<m^*+{1\over k}, \quad k\in \mathbb{N}^*.$$
Since $\{ \zeta _k\} _{k\in \mathbb{N}^*}\subset {\cal M}$ it follows that there exists a subsequence (also denoted by $\{ \zeta_k\}$) such that
$$\zeta _k \longrightarrow \zeta ^*\in {\cal M} \ \mbox{\rm  weakly * in } L^{\infty}(\mathbb{R}^d),\quad {\cal K}(\zeta _k)\longrightarrow {\cal K}(\zeta ^*),  \ \mbox{\rm uniformly for } x\in \mathbb{R}^d.$$
Using the convexity of $Q$ with respect to $\zeta $, we obtain that
$$I(\zeta _k)\geq \int_0^T\int_{\mathbb{R}^d}G(t,x)\rho ^{\zeta _k}(t,x)dx \ dt+\int_{\mathbb{R}^d}G_T(x)\rho ^{\zeta _k}(T,x)dx$$
$$+\int_{\mathbb{R}^d}Q(x,\zeta ^*(x))dx+\int_{\mathbb{R}^d}{{\partial Q}\over {\partial z}}(x,\zeta ^*(x))(\zeta _k(x)-\zeta ^*(x))dx .$$
The last inequality implies via Lemma 3.2 that
$$m^*\geq \liminf _{k\rightarrow +\infty }I(\zeta _k)\geq I(\zeta^*),$$
and so $\zeta ^*$ is an optimal control for problem $(P)$.
\vspace{3mm}

\noindent
\it Proof of Lemma 3.1. \rm 
Consider for the beginning the function $\Phi $ defined in the proof of Lemma 3.3 in \cite{BR4}. Actually, we recall that this particular function satisfies
$\Phi \in C^2(\mathbb{R}^d)$, $1\leq \Phi (x)$, $\forall x\in \mathbb{R}^d$ and 
$$\lim_{|x|_d\rightarrow +\infty }\Phi (x)=+\infty, \quad \nabla \Phi \in L^{\infty}(\mathbb{R}^d;\mathbb{R}^d), \quad \Delta \Phi \in L^{\infty}(\mathbb{R}^d).$$

\begin{itemize}
\item[Case 1.] Let $f \in L^1(\mathbb{R}^d )\cap L^{\infty}(\mathbb{R}^d)$ such that $\int_{\mathbb{R}^d}|f(x)|\Phi (x) dx < +\infty$.

From the proof of Lemma 2.5 in \cite{BR5} we have that $y_n$, the unique solution to (\ref{eq3.2}), satisfies $\beta (y_n) \in H^1(\mathbb{R}^d)$.

By equation (\ref{eq3.2}) we get that
$$\int_{\mathbb{R}^d}y_n\beta(y_n) dx +h\int_{\mathbb{R}^d}|\nabla \beta (y_n)|_d^2 dx-h\int_{\mathbb{R}^d}{\cal K}(\zeta _n)b(y_n)y_n\cdot \nabla \beta(y_n) dx =\int_{\mathbb{R}^d}f \beta(y_n) dx.$$
Notice that $\beta(y_n)=\beta(y_n)-\beta(0)=\beta '(\xi_n)y_n$ and $\nabla \beta(y_n)=\beta '(y_n)\nabla y_n$ (where $\xi _n(x)$ is an intermediate point).
It follows that for $h>0$ sufficiently small, $\{\beta(y_n)\}_{n\in \mathbb{N}^*}$ is bounded in $H^1(\mathbb{R}^d)$ and $\{y_n\}_{n\in \mathbb{N}^*}$ is bounded in $H^1(\mathbb{R}^d)$ as well.

\noindent
It follows that there exists a subsequence $y_{n_k}\longrightarrow \tilde{y}$ in $L^2(B(0_d;1))$, there exists a subsequence $y_{n_{k_l}}\longrightarrow \tilde{y}$ in $L^2(B(0_d;2))$, ...

Taking the diagonal subsequence (also denoted by $\{y_n\}_{n\in \mathbb{N}^*}$) we get that 
$y_n \longrightarrow \tilde{y}$ in $L^2_{loc}(\mathbb{R}^d)$.
Repeating the argument in Lemma 3.3 in \cite{BR4} we obtain that there exists a nonnegative constant $C$ such that
$$\int_{\mathbb{R}^d}|y_n(x)|\Phi (x) dx \leq C, \quad \forall n\in \mathbb{N}^*.$$
Using Fatou's lemma we obtain $\displaystyle \int_{\mathbb{R}^d}|\tilde{y}(x)|\Phi (x) dx \leq C$.

Since $\displaystyle \int_{\mathbb{R}^d}[y_n\varphi-h\beta(y_n)\Delta \varphi-h{\cal K}(\zeta _n)b^*(y_n)\cdot \nabla \varphi] dx=\int_{\mathbb{R}^d}f\varphi dx$, \quad $\forall \varphi \in C^{\infty}_0(\mathbb{R}^d)$,
we get that
$$\int_{\mathbb{R}^d}[\tilde{y}\varphi-h\beta(\tilde{y})\Delta \varphi-h{\cal K}(\zeta )b^*(\tilde{y})\cdot \nabla \varphi] dx=\int_{\mathbb{R}^d}f\varphi dx.$$
So, $\tilde{y}$ satisfies (\ref{eq3.3}) in distributional sense.

On the other hand, by Lemma 2.5 in \cite{BR5} we have that
$\displaystyle \int_{\mathbb{R}^d}|y_n|dx\leq \int_{\mathbb{R}^d}|f|dx$ and consequently $\displaystyle \int_{\mathbb{R}^d}|\tilde{y}|dx\leq \int_{\mathbb{R}^d}|f|dx$.

We may infer that
$\tilde{y} \in L^1 (\mathbb{R}^d), \ \beta(\tilde{y}) \in L^1(\mathbb{R}^d)$ (and so $-\Delta \beta(\tilde{y})+\nabla \cdot ({\cal K}(\zeta )b^*(\tilde{y})) \in L^1(\mathbb{R}^d)$).
So, $\tilde{y}$ is solution to (\ref{eq3.3}) and by Lemma 2.5 in \cite{BR5} we conclude that $\tilde{y}=y$.
\vspace{2mm}

Let $\tilde{\Phi}_{\cal R}=\min\{ \Phi(x); |x|_d \geq {\cal R}\}$, ${\cal R}>0$. This yields 
$$\int_{\mathbb{R}^d}|y_n-y| dx=\int_{B(0_d;{\cal R})}|y_n-y| dx+\int_{|x|_d\geq {\cal R}}|y_n-y| dx.$$
Since 
$$\displaystyle \int_{|x|_d\geq {\cal R}}|y_n|\tilde{\Phi}_{\cal R} dx\leq \int_{|x|_d\geq {\cal R}}|y_n| \Phi (x) dx \leq  \int_{\mathbb{R}^d}|y_n| \Phi (x) dx \leq C,$$
it follows that
$$ \int_{|x|_d\geq {\cal R}}|y_n| dx \leq {{C}\over {\tilde{\Phi}_{\cal R}}} \ \mbox{\rm and so } \int_{|x|_d\geq {\cal R}}|y| dx \leq {{C}\over {\tilde{\Phi}_{\cal R}}}.$$

\noindent
On the other hand, since $\tilde{\Phi}_{\cal R} \rightarrow +\infty$ as ${\cal R} \rightarrow +\infty$, it follows that for any $\varepsilon >0$, there exists for ${\cal R}>0$ sufficiently large such that
$$\int_{|x|_d\geq {\cal R}}|y_n-y| dx \leq \int_{|x|_d\geq {\cal R}}(|y_n|+|y|) dx \leq 2\cdot {{\varepsilon}\over 2}=\varepsilon,$$
and so
$$\limsup _{n\rightarrow +\infty }\int_{\mathbb{R}^d}|y_n-y| dx \leq \varepsilon,  \ \forall \varepsilon >0, \quad \mbox{\rm and  } y_n \longrightarrow y \quad \mbox{\rm in } L^1(\mathbb{R}^d).$$
\vspace{1mm}

\noindent
\bf Remark 3.2. \rm The conclusion holds not only for a subsequence of $\{ y_n\}_{n\in \mathbb{N}^*}$, but for the sequence itself. We argue by contradiction.

Assume that there exists a subsequence $\{y_{n_l}\}_{l\in \mathbb{N}^*}$ such that $\|y_{n_l}-y\|_{L^1(\mathbb{R}^d)}\rightarrow \theta _0>0$.
The first part of the proof allows us to conclude that on a subsequence $\{y_{n_{l_q}}\}_{q\in \mathbb{N}^*}$ we have that $y_{n_{l_q}} \longrightarrow y$ in ${L^1(\mathbb{R}^d)}$ and so $ \|y_{n_{l_q}}-y\|_{L^1(\mathbb{R}^d)}\longrightarrow 0$; absurd.
\vspace{2mm}

\item[Case 2.] Let $f \in L^1(\mathbb{R}^d)$ and $\{f_k\}_{k\in \mathbb{N}^*}$ satisfying the hypotheses of $f$ in case 1 and satisfying $f_k \longrightarrow f$ in $L^1(\mathbb{R}^d)$. Let
\begin{itemize}
\item[$\bullet $] $y^k_n$ solution to (\ref{eq3.3}) with $\zeta :=\zeta _n$, $f:=f_k$;
\item[$\bullet $] $y_n$ solution to (\ref{eq3.3}) with $\zeta :=\zeta _n$;
\item[$\bullet $] $y^k$ solution to (\ref{eq3.3}) with $f:=f_k$.
\end{itemize}
We have that
$$\|y_n-y\|_{L^1(\mathbb{R}^d)} \leq \|y_n-y^k_n\|_{L^1(\mathbb{R}^d)}+\|y^k_n-y^k\|_{L^1(\mathbb{R}^d)}+\|y^k-y\|_{L^1(\mathbb{R}^d)}.$$
Since
$$\|y_n-y^k_n\|_{L^1(\mathbb{R}^d)} \leq \|f-f_k\|_{L^1(\mathbb{R}^d)}, \quad \|y^k-y\|_{L^1(\mathbb{R}^d)} \leq \|f-f_k\|_{L^1(\mathbb{R}^d)},$$
we conclude that
$$\|y_n-y\|_{L^1(\mathbb{R}^d)} \leq 2\|f-f_k\|_{L^1(\mathbb{R}^d)}+\|y^k_n-y^k\|_{L^1(\mathbb{R}^d)}.$$
For a fixed $k$ we have that $y^k_n-y^k \longrightarrow 0$ in $L^1(\mathbb{R}^d)$.

Hence $\limsup_{n\rightarrow +\infty } \|y_n-y\|_{L^1(\mathbb{R}^d)} \leq 2 \|f-f_k\|_{L^1(\mathbb{R}^d)}$, $\forall k\in \mathbb{N}^*$, and so
$$\limsup_{n\rightarrow +\infty } \|y_n-y\|_{L^1(\mathbb{R}^d)} \leq 0 \ \mbox{\rm and } y_n \longrightarrow y \ \mbox{\rm in } L^1(\mathbb{R}^d).$$
\end{itemize}

\section{A penalized optimal control problem ${\bf (P_h)}$}

For any $N\in \mathbb{N}^*$ we consider $h={T\over N}$ and the following optimal control problem:
$$\underset{\zeta \in {\cal M}}{\text{Minimize}} \  \int_0^T \int_{\mathbb{R}^d}G(t,x)\rho ^{\zeta }_h(t,x) dx \ dt+\int_{\mathbb{R}^d}G_T(x)\rho ^{\zeta }_h(T,x)dx \leqno({\bf P_h})$$
$$+\int_{\mathbb{R}^d}Q(x,\zeta (x))dx+{1\over 2}\int_{\mathbb{R}^d}|\zeta (x)-\zeta ^*(x)|^2dx, $$
where $\rho ^{\zeta }_h$ and $\rho ^{\zeta ,i}_h$ were defined in section 2, and $\zeta ^*$ is an optimal control for problem $(P)$.
\vspace{2mm}

\noindent
\bf Remark 4.1. \rm The following equality holds
$$\int_0^T \int_{\mathbb{R}^d}G(t,x)\rho ^{\zeta }_h(t,x) dx \ dt+\int_{\mathbb{R}^d}G_T(x)\rho ^{\zeta }_h(T,x)dx$$
$$=\sum_{i=1}^N \int_{\mathbb{R}^d}\tilde{G}^i(x)\rho ^{\zeta ,i}_h(x)dx+\int_{\mathbb{R}^d} G_T (x)\rho ^{\zeta ,N}_h(x)dx,$$
where
$$\tilde{G}^i(x)=\int _{(i-1)h}^{ih}G(t,x) dt, \quad  i=1,2,...,N.$$
\vspace{4mm}

For the remainder of the section we consider $h={T\over N}$ to be sufficiently small.

\subsection{Existence of an optimal control for $\bf (P_h)$}

\noindent
\bf Lemma 4.1. \it If $\{ \zeta_n \}_{n\in \mathbb{N}^*} \subset {\cal M}$, $\zeta _n \longrightarrow \zeta \in {\cal M}$ weakly * in $L^{\infty}(\mathbb{R}^d)$,
${\cal K}(\zeta _n)\longrightarrow {\cal K}(\zeta )$, uniformly for $x\in \mathbb{R}^d,$
and if $\{f_n\}_{n\in \mathbb{N}^*}\subset L^1(\mathbb{R}^d)$, $f_n \longrightarrow f$ in $L^1(\mathbb{R}^d)$, then
$$y^n_n \longrightarrow y \quad \mbox{\it in } L^1(\mathbb{R}^d),$$
where $y_n^k, \ y_n$ where introduced in section 3 (in the proof of Lemma 3.1).
\vspace{2mm}

\noindent
Proof. \rm By Lemma 2.5 in \cite{BR5} we get that
$$\|y^n_n-y\|_{L^1(\mathbb{R}^d)} \leq \|y^n_n-y_n\|_{L^1(\mathbb{R}^d)}+\|y_n-y\|_{L^1(\mathbb{R}^d)} \leq \|f_n-f\|_{L^1(\mathbb{R}^d)}+\|y_n-y\|_{L^1(\mathbb{R}^d)}.$$
Using the convergence of $\{ f_n\}_{n\in \mathbb{N}^*}$ and that $y_n\longrightarrow y$ in $L^1(\mathbb{R}^d)$ (by Lemma 3.1) we get the conclusion. 
\vspace{3mm}

\noindent
\bf Lemma 4.2. \it If $\{ \zeta_n \}_{n\in \mathbb{N}^*} \subset {\cal M}$, $\zeta _n \longrightarrow \zeta \in {\cal M}$ weakly * in $L^{\infty}(\mathbb{R}^d)$,
${\cal K}(\zeta _n)\longrightarrow {\cal K}(\zeta )$, uniformly for $x\in \mathbb{R}^d,$
and if 
$(\rho ^{\zeta _n,i}_h)_{i=0,1,...,N}$ and $(\rho ^{\zeta ,i}_h)_{i=0,1,...,N}$ are defined by $(\ref{eq2.7})$, then
$$\rho ^{\zeta _n,i}_h \longrightarrow \rho ^{\zeta ,i}_h \quad \mbox{\it in } L^1(\mathbb{R}^d), \quad \forall i \in \{0,1,...,N\}.$$

\noindent
Proof. \rm This follows by finite induction and Lemma 4.1.
\vspace{3mm}

\noindent
\bf Theorem 4.3. \it Problem $(P_h)$ has at least one optimal control $\zeta ^*_h$. 
\vspace{2mm}

\noindent
Proof. \rm Denote by $I_h$ the cost functional associated to problem $(P_h)$. 
Let $m^*_h = \underset{\zeta \in {\cal M}}{\text{inf}}I_h(\zeta ) \in [0,+\infty ), $ and let
$\{ \zeta_n \}_{n\in \mathbb{N}^*} \subset {\cal M}$, $\zeta _n \longrightarrow \zeta ^*_h\in {\cal M}$ weakly * in $L^{\infty}(\mathbb{R}^d)$,
$${\cal K}(\zeta _n)\longrightarrow {\cal K}(\zeta _h^*),  \quad \mbox{\rm uniformly for } x\in \mathbb{R}^d,$$
such that:
$$m^*_h \leq I_h(\zeta _n) < m^*_h +{1\over n}, \quad n\in \mathbb{N}^*.$$
By Lemma 4.2 we have that 
$$\int_{\mathbb{R}^d}G_T \rho ^{\zeta _n,N}_h dx \longrightarrow \int_{\mathbb{R}^d}G_T \rho ^{\zeta ^*_h,N}_h dx,$$
and 
$$\sum_{i=1}^N \int_{\mathbb{R}^d}\tilde{G}^i \rho ^{\zeta _n,i}_h dx \longrightarrow \sum_{i=1}^N \int_{\mathbb{R}^d}\tilde{G}^i \rho ^{\zeta ^*_h,i}_h dx.$$
On the other hand (from the convexity of $Q$ with respect to $\zeta$)
$$\int_{\mathbb{R}^d}Q(x,\zeta _n (x))dx+{1\over 2}\int_{\mathbb{R}^d}|\zeta _n(x)-\zeta ^*(x)|^2dx \geq \int_{\mathbb{R}^d}Q(x,\zeta _h^* (x))dx$$
$$+\int_{\mathbb{R}^d}{{\partial Q}\over{\partial z}}(x,\zeta _h^*(x))(\zeta _n(x)-\zeta _h^*(x))dx.$$
This implies that
$$m^*_h\geq \lim _{n\rightarrow +\infty } I_h(\zeta _n)=\liminf _{n\rightarrow +\infty } I_h(\zeta _n)\geq I(\zeta _h^*)\geq m_h^*.$$
We conclude that $\zeta ^*_h$ is an optimal control for problem $(P_h)$.

\subsection{Necessary optimality conditions for $\bf (P_h)$}
\large 

Let $\zeta _h^*$ be an optimal control for problem $(P_h)$. 
Let $\xi \in L^{\infty }(\mathbb{R}^d)$ be such that
$\zeta_h^*+\varepsilon \xi\in {\cal M}$, for any $\varepsilon >0$ sufficiently small.
\vspace{3mm}

For any $\zeta \in {\cal M}$, let $z^{\zeta ,i+1}$ be the weak solution to
\begin{equation}\label{ec_zh}
\left\{ \begin{array}{ll}
z^{\zeta ,i+1}-h\Delta (\beta '(\rho _h^{\zeta ,i+1})z^{\zeta ,i+1})+h\nabla \cdot ({\cal K}(\zeta )(b^*)'(\rho _h^{\zeta ,i+1})z^{\zeta ,i+1})+h\nabla \cdot ({\cal K}(\xi )b^*(\rho _h^{\zeta .i+1}))=z^{\zeta ,i}, \\
\hspace{10cm} i=0,1,...,N-1, \\
z^{\zeta ,0}=0.
          \end{array}
\right.
\end{equation}
Actually, $z^{\zeta ,i+1}$ is a weak solution to the equation in (\ref{ec_zh}) if $\tilde{z}^{i+1}\in H^1(\mathbb{R}^d)$ is the unique weak solution to
$${1\over {\beta '(\rho_h^{\zeta ,i+1})}}\tilde{z}^{i+1}-h\Delta \tilde{z}^{i+1}+h\nabla \cdot \left( {\cal K}(\zeta ){{(b^*)'(\rho _h^{\zeta ,i+1})}\over {\beta '(\rho _h^{\zeta ,i+1})}}\tilde{z}^{i+1}\right) =-h\nabla \cdot ({\cal K}(\xi )b^*(\rho _h^{\zeta .i+1}))+{1\over {\beta '(\rho_h^{\zeta ,i})}}\tilde{z}^{i} $$
and $\tilde{z}^0=0$. Here $\tilde{z}^i=\beta '(\rho _h^{\zeta ,i})z^{\zeta ,i}$.
The existence and uniqueness of a weak solution to this equation follows via Lax-Milgram's lemma (for sufficiently small $h$).
\vspace{2mm}

We postpone the proof of the following lemma.
\vspace{3mm}

\noindent
\bf Lemma 4.4. \it The following convergences hold for any $i\in \{ 1,2,...,N\}$:
\begin{itemize}
\item[(i)] $\rho _h^{\zeta ^*_h+\varepsilon \xi ,i}\longrightarrow \rho _h^{\zeta ^*_h,i} \quad \mbox{\it in } L^2(\mathbb{R}^d), \ \mbox{\it as } \varepsilon\longrightarrow 0+$;
\item[(ii)] ${{\rho _h^{\zeta ^*_h+\varepsilon \xi ,i}-\rho _h^{\zeta ^*_h,i}}\over {\varepsilon }}\longrightarrow z^{\zeta ^*_h,i} \quad \mbox{\it in } L^2(\mathbb{R}^d), \ \mbox{\it as } \varepsilon\longrightarrow 0+$
(on a subsequence).
\end{itemize}
\rm
\vspace{3mm}

\noindent
\bf Theorem 4.5. \it If $p_h^{\zeta ^*_h,i+1}$ is the unique weak solution to
\begin{equation}\label{eq4.1}
\left\{ \begin{array}{ll}
p^{\zeta ,i}-h\beta '(\rho _h^{\zeta ,i+1})\Delta p^{\zeta ,i}-h{\cal K}(\zeta )(b^*)'(\rho _h^{\zeta ,i+1})\cdot \nabla p^{\zeta ,i}+\tilde{G}^{i+1}=p^{\zeta ,i+1}, \\
\hspace{10cm} i=N-1,...,1,0, \\
p^{\zeta ,N}=-G_T,
          \end{array}
\right.
\end{equation}
corresponding to $\zeta :=\zeta ^*_h$, then
\begin{equation}\label{eq4.2}
\int_{\mathbb{R}^d}{\cal K}(\xi )\cdot [\sum_{i=1}^N[hb^*(\rho _h^{\zeta ^*_h,i})\nabla p_h^{\zeta ^*_h,i-1}]dx
-\int_{\mathbb{R}^d}\xi [{{\partial Q}\over {\partial z}}(x,\zeta _h^*(x))+(\zeta ^*_h(x)-\zeta ^*(x))]dx\leq 0,
\end{equation}
for any $\xi \in L^{\infty }(\mathbb{R}^d)$ such that
$\zeta_h^*+\varepsilon \xi\in {\cal M}$, for any $\varepsilon >0$ sufficiently small.
\rm
\vspace{6mm}

We say that $p^{\zeta ,i}\in H^1(\mathbb{R}^d)$ is a weak solution to the equation in (\ref{eq4.1}) if it is a weak solution to 
$${1\over {\beta '(\rho _h^{\zeta ,i+1})}}p^{\zeta ,i}-h\Delta p^{\zeta ,i}-h{\cal K}(\zeta ){{(b^*)'(\rho _h^{\zeta ,i+1})}\over {\beta '(\rho _h^{\zeta ,i+1})}}\cdot \nabla p^{\zeta ,i}+{1\over {\beta '(\rho _h^{\zeta ,i+1})}}\tilde{G}^{i+1}={1\over {\beta '(\rho _h^{\zeta ,i+1})}}p^{\zeta ,i+1}.$$
The existence and uniqueness of a weak solution follows via Lax-Milgram's lemma (for sufficiently small $h$).
\vspace{2mm}

\noindent
\bf Remark 4.2. \rm We have that $\zeta _h^*(x)=0$ a.e. $x\in Ext (B(0_d;R_0))$. On the other hand for almost any $x\in B(0_d;R_0)$ we have that
$$\zeta _h^*(x)=\left\{ \begin{array}{ll}0, & \mbox{\rm if } \displaystyle \int_{\mathbb{R}^d}\nabla {\cal G}_R(x'-x)\cdot \sum_{i=1}^N[hb^*(\rho _h^{\zeta ^*_h,i})\nabla p_h^{\zeta ^*_h,i-1}](x')dx' \\
~ & \displaystyle \ \ \ \ \ +{{\partial Q}\over {\partial z}}(x,\zeta _h^*(x))+(\zeta ^*_h(x)-\zeta ^*(x)) >0 \\
~  & ~ \\
\tilde{M}_0,  & \mbox{\rm if } \displaystyle \int_{\mathbb{R}^d}\nabla {\cal G}_R(x'-x)\cdot \sum_{i=1}^N[hb^*(\rho _h^{\zeta ^*_h,i})\nabla p_h^{\zeta ^*_h,i-1}](x')dx' \\
~ & \ \ \ \ \ \displaystyle +{{\partial Q}\over {\partial z}}(x,\zeta _h^*(x))+(\zeta ^*_h(x)-\zeta ^*(x)) <0 .
\end{array}
\right. $$
\vspace{2mm}

\noindent
\it Proof of Theorem 4.5. \rm Since $\zeta ^*_h$ is an optimal control for problem $(P_h)$, then $I_h(\zeta ^*_h)\leq I_h(\zeta ^*_h+\varepsilon \xi )$, for any $\varepsilon >0$ sufficiently small.
After an easy evaluation and using the Lemma 4.4 we get that
\begin{equation}\label{ineq}
0\leq \sum_{i=1}^N\int_{\mathbb{R}^d}\tilde{G}^{i}z^{\zeta ^*_h,i}dx+\int_{\mathbb{R}^d}G_Tz^{\zeta ^*_h,N}dx+\int_{\mathbb{R}^d}[{{\partial Q}\over {\partial z}}(x,\zeta _h^*(x))+(\zeta _h^*(x)-\zeta ^*(x))] \xi (x)dx.
\end{equation}
We multiply (\ref{eq4.1}) by $z^{\zeta ^*_h,i+1}$ and integrate over $\mathbb{R}^d$. After some evaluation and taking into account (\ref{ec_zh}) we obtain that
$$\int_{\mathbb{R}^d}\tilde{G}^{i+1}z^{\zeta ^*_h,i+1}dx=\int_{\mathbb{R}^d}p_h^{\zeta ^*_h,i+1}z^{\zeta ^*_h,i+1}dx+\int_{\mathbb{R}^d}p_h^{\zeta ^*_h,i}[-z^{\zeta ^*_h,i}+h\nabla \cdot ({\cal K}(\xi )b^*(\rho _h^{\zeta ^*_h,i+1}))]dx.$$
By summation and using that $z^{\zeta ^*_h,0}=0$ we obtain that
$$\sum_{i=1}^N \int_{\mathbb{R}^d}\tilde{G}^{i}z^{\zeta ^*_h,i}dx=-\int_{\mathbb{R}^d}G_Tz^{\zeta ^*_h,N}dx-\sum_{i=1}^N\int_{\mathbb{R}^d}hb^*(\rho _h^{\zeta ^*_h,i}){\cal K}(\xi )\cdot \nabla p_h^{\zeta ^*_h,i-1}dx$$
and using (\ref{ineq}) we get (\ref{eq4.2}).
\vspace{3mm}

\noindent
\it Proof of Lemma 4.4. \rm  By Lemma 3.1 in \cite{BR4} it follows that there exists $C>0$ such that
$$\| \rho_h^{\zeta ^*_h+\varepsilon \xi ,i}\| _{L^{\infty }(\mathbb{R}^d)}\leq C, \quad \forall 0\leq \varepsilon \leq \varepsilon _0, \ i=0,1,...,N.$$
Since $\rho _h^{\zeta ^*_h+\varepsilon \xi,i}\longrightarrow \rho _h^{\zeta ^*_h,i} \ \mbox{\rm in } L^1(\mathbb{R}^d), \ i=0,1,...,N$, as $\varepsilon \longrightarrow 0+$  (by Lemma 4.2) we conclude that
(i) is satisfied.
\vspace{2mm}

Let $y^{\varepsilon ,i}=\beta (\rho _h^{\zeta ^*_h+\varepsilon \xi ,i}) , \quad y^{i}=\beta (\rho _h^{\zeta ^*_h,i})$, 
$\eta ^i=\beta '(\rho _h^{\zeta ^*_h,i})z^{\zeta ^*_h,i}$, $\eta _{\varepsilon }^i=y^{\varepsilon ,i}-y^i$, $\tilde{\eta }_{\varepsilon }^i={{\eta _{\varepsilon }^i}\over {\varepsilon }}$.

\noindent
We have that $\zeta _{\varepsilon }^{i+1}$ is the weak solution to
$$\beta ^{-1}(y^{\varepsilon ,i+1})-\beta ^{-1}(y^{i+1})-h\Delta \eta _{\varepsilon }^{i+1}+h\nabla \cdot ({\cal K}(\zeta _h^*)(b^*(\beta ^{-1}(y^{\varepsilon ,i+1}))-b^*(\beta ^{-1}(y^{i+1}))))$$
$$+\varepsilon h\nabla \cdot ({\cal K}(\xi )b^*(\beta ^{-1}(y^{\varepsilon ,i+1})))=\beta ^{-1}(y^{\varepsilon ,i})-\beta ^{-1}(y^{i}).$$
If we denote by $w_{\varepsilon }^{i}=\tilde{\eta }_{\varepsilon }^i-\eta ^{i}$, then we have that $w_{\varepsilon }^{i+1}$ is the unique weak solution to
$${{\beta ^{-1}(y^{\varepsilon ,i+1})-\beta ^{-1}(y^{i+1})}\over {\varepsilon }}-{1\over {\beta '(\beta ^{-1}(y^{i+1}))}}\eta ^{i+1}-h\Delta w_{\varepsilon }^{i+1}$$
$$+h\nabla \cdot \left({\cal K}(\zeta _h^*) ({{b^*(\beta ^{-1}(y^{\varepsilon ,i+1}))-b^*(\beta ^{-1}(y^{i+1}))}\over {\varepsilon }}-{{(b^*)'(\beta ^{-1}(y^{i+1}))}\over {\beta '(\beta ^{-1}(y^{i+1}))}}\eta ^{i+1})\right) $$
$$={{\beta ^{-1}(y^{\varepsilon ,i})-\beta ^{-1}(y^{i})}\over {\varepsilon }}-{1\over {\beta '(\beta ^{-1}(y^{i}))}}\eta ^{i}.$$
If we apply the Lagrange's mean value theorem we get after some evaluation and by finite induction that for a subsequence (also indexed by $\varepsilon $):
$$w_{\varepsilon }^i\longrightarrow 0 \quad \mbox{\rm in } H^1(\mathbb{R}^d), \quad \forall i=0,1,...,N,$$
and so
$${{\beta (\rho _h^{\zeta ^*_h+\varepsilon \xi ,i})-\beta (\rho _h^{\zeta ^*_h,i})}\over {\varepsilon }}\longrightarrow \beta '(\rho _h^{\zeta ^*_h,i})z^{\zeta ^*_h,i} \quad \mbox{\rm in } L^2(\mathbb{R}^d), \quad \mbox{\rm as } \varepsilon \longrightarrow 0+, \ \forall i=0,1,...,N.$$
If we apply once more the Lagrange's mean value theorem we get that on a subsequence (ii) holds.
\vspace{3mm}

\noindent
\bf Remark 4.3. \rm Note that for any $\zeta \in {\cal M}$, (\ref{eq2.7}) is a backward Euler scheme for (\ref{nfp1}). Moreover, for any $\zeta \in {\cal M}$ we get that
$$\rho _h^{\zeta }(t)\longrightarrow \rho ^{\zeta }(t) \quad \mbox{\rm in } L^1(\mathbb{R}^d), \ \mbox{\rm as } h\longrightarrow 0,$$
uniformly for $t\in [0,T]$.

\section{Relationship between problems $\bf (P)$ and $\bf (P_h)$. Optimality conditions for $\bf (P)$}

For any $h=\displaystyle {T\over N}$, $N\in \mathbb{N}^*$, we consider $\zeta ^*_h\in {\cal M}$ an optimal control for problem $(P_h)$. 
\vspace{3mm}

\noindent
\bf Lemma 5.1. \it The sequence
$\{ \rho _h^{\zeta }\}$ is bounded in $L^{\infty }(0,T;L^2(\mathbb{R}^d))$, uniformly for $\zeta \in {\cal M}$ and
$\{\nabla \rho _h^{\zeta }\} $  is bounded in $L^2(0,T;L^2(\mathbb{R}^d)^d)$, uniformly for $\zeta \in {\cal M}$, for $h>0$ sufficiently small.
\vspace{2mm}

\noindent
Proof. \rm By (\ref{eq2.7}) we obtain that for any $\zeta \in {\cal M}$:
$$\int_{\mathbb{R}^d}|\rho _h^{\zeta ,i+1}|^2dx +h\int_{\mathbb{R}^d}\nabla \beta (\rho _h^{\zeta ,i+1})\cdot \nabla \rho _h^{\zeta ,i+1}dx=\int_{\mathbb{R}^d}\rho_h^{\zeta ,i}\rho _h^{\zeta ,i+1}dx+h\int_{\mathbb{R}^d}{\cal K}(\zeta )b^*(\rho _h^{\zeta ,i+1})\cdot \nabla \rho _h^{\zeta ,i+1}dx,$$
and so
$${1\over 2}\int_{\mathbb{R}^d}|\rho _h^{\zeta ,i+1}|^2dx +h\gamma _0\int_{\mathbb{R}^d}|\nabla \rho _h^{\zeta ,i+1}|_d^2dx\leq
{1\over 2}\int_{\mathbb{R}^d}|\rho_h^{\zeta ,i}|^2dx+h\tilde{M}_1\int_{\mathbb{R}^d}|\rho _h^{\zeta ,i+1}| |\nabla \rho _h^{\zeta ,i+1}|_ddx$$
where $\tilde{M}_1$ is a nonnegative constant. We immediately obtain that there exists a positive constant $M_2$ such that
$${1\over 2}\int_{\mathbb{R}^d}|\rho _h^{\zeta ,i+1}|^2dx +{{h\gamma _0}\over 2}\int_{\mathbb{R}^d}|\nabla \rho _h^{\zeta ,i+1}|_d^2dx\leq
{1\over 2}\int_{\mathbb{R}^d}|\rho_h^{\zeta ,i}|^2dx+{{hM_2}\over 2}\int_{\mathbb{R}^d}|\rho _h^{\zeta ,i+1}|^2dx,$$
and consequently for $h>0$ sufficiently small
$$(1-hM_2)\int_{\mathbb{R}^d}|\rho _h^{\zeta ,i+1}|^2dx+h\gamma _0\int_{\mathbb{R}^d}|\nabla \rho _h^{\zeta ,i+1}|_d^2dx\leq\int_{\mathbb{R}^d}|\rho_h^{\zeta ,i}|^2dx.$$
After an easy evaluation we get that
$$(1-hM_2)^N\| \rho _h^{\zeta }(t)\|^2_{L^2}+\gamma _0(1-hM_2)^N\int_0^t\int_{\mathbb{R}^d}|\nabla \rho _h^{\zeta }(s,x)|^2_ddx \ ds\leq \|\rho _0\| ^2_{L^2},$$
for any $t\in [0,T], \ \zeta \in {\cal M}.$
Since $h={T\over N}$ we get the conclusion of the lemma.
\vspace{3mm}

We postpone the proof ot the following result.
\vspace{3mm}

\noindent
\bf Lemma 5.2. \it If for a sequence $\{ \zeta  _h\}_{h>0}\subset {\cal M}$ we have that 
$\zeta _{h}\longrightarrow \zeta \in {\cal M}$ weakly * in $L^{\infty }(\mathbb{R}^d),$
$${\cal K}(\zeta _h)\longrightarrow {\cal K}(\zeta ) \ \mbox{\it uniformly in } \mathbb{R}^d,\quad \nabla {\cal K}(\zeta _h)(x)\longrightarrow \nabla {\cal K}(\zeta ) (x), \ \forall x \in \mathbb{R}^d,$$
as $h\rightarrow 0$, then
$$\rho _h^{\zeta _h}(t)\longrightarrow \rho ^{\zeta }(t) \ in \ L^1(\mathbb{R}^d), \ \mbox{\it uniformly for } t\in [0,T].$$
\vspace{3mm}

\noindent
\bf Theorem 5.3 (The relationship between $\bf (P_h)$ and $\bf (P)$). \it The following convergences hold
\begin{itemize}
\item[(j)] $\zeta_h^*\longrightarrow \zeta ^* $ in $L^2(\mathbb{R}^d)$;
\item[(jj)] $I_h(\zeta _h^*)\longrightarrow I(\zeta ^*)$;
\item[(jjj)] $I(\zeta _h^*)\longrightarrow I(\zeta ^*)$, as $h\rightarrow 0$.
\end{itemize}
\vspace{1mm}
\rm

\noindent
Proof. \rm Since $\{ \zeta ^*_h\} \subset {\cal M} $ it follows that there exists a subsequence (also denoted by $\{ \zeta ^*_h\}$)
such that $\zeta ^*_h\longrightarrow \tilde{\zeta }\in {\cal M}$ weakly * in $L^{\infty }(\mathbb{R}^d)$ and
$$\rho _h^{\zeta ^*_h}(t)\longrightarrow \rho^{\tilde{\zeta }}(t) \quad \mbox{\rm in } L^1(\mathbb{R}^d), \ \mbox{\rm uniformly for } t\in [0,T].$$
This yields
$$I_h(\zeta ^*)\geq I_h(\zeta ^*_h)\geq \int_0^T \int_{\mathbb{R}^d}G(t,x)\rho _h^{\zeta ^*_h}(t,x) dx \ dt+\int_{\mathbb{R}^d}G_T(x)\rho _h^{\zeta ^*_h}(T,x)dx$$
$$+\int_{\mathbb{R}^d}Q(x,\tilde{\zeta }(x))dx+\int_{\mathbb{R}^d}{{\partial Q}\over {\partial z}}(x,\tilde{\zeta }(x))(\zeta ^*_h(x)-\tilde{\zeta }(x))dx+{1\over 2}\int_{\mathbb{R}^d}|\zeta ^*_h(x)-\zeta ^*(x)|_d^2dx .$$
By Crandall-Liggett's theorem we have that $\rho _h^{\zeta ^*}(t)\longrightarrow \rho ^{\zeta ^*}(t)$ in $L^1(\mathbb{R}^d)$, uniformly for $t\in [0,T]$, and so we get that
$$I_h(\zeta ^*)\longrightarrow I(\zeta ^*),$$
and
$$I(\zeta ^*)\geq \limsup _{h\rightarrow 0}I_h(\zeta ^*_h)\geq \liminf _{h\rightarrow 0}I_h(\zeta ^*_h)\geq I(\tilde{\zeta})+{1\over 2}\int_{\mathbb{R}^d}|\tilde{\zeta }(x)-\zeta ^*(x)|^2_ddx \geq I(\zeta ^*).$$
Since $\zeta ^*$ is an optimal control for problem $(P)$ we may infer that $\tilde{\zeta }(x)=\zeta ^*(x)$ a.e. $x\in \mathbb{R}^d$ and that on a subsequence (j) and (jj) hold.

\noindent
Arguing by contradiction we get that (j) and (jj) hold in general (not only on a subsequence).
\vspace{2mm}

On the other hand, using that $\rho ^{\zeta ^*_h}(t)\longrightarrow \rho ^{\zeta ^*}(t)$ in $L^1(\mathbb{R}^d)$, uniformly for $t\in [0,T]$ (by Lemma 3.2), we get that (jjj) holds.
\vspace{3mm}

\noindent
\bf Theorem 5.4 (Necessary optimality conditions for $\bf (P)$). \it If $p\in W^{1,2}([0,T];L^2(\mathbb{R}^d))\cap L^2(0,T;H^2(\mathbb{R}^d))$ satisfies
\begin{equation}\label{ec5.1}
\left\{ \begin{array}{ll}
\displaystyle
{{dp}\over {dt}}(t)=-\beta '(\rho ^{\zeta ^*}(t))\Delta p(t)-{\cal K}(\zeta ^*)(b^*)'(\rho ^{\zeta ^*}(t))\cdot \nabla p(t)+G(t), \quad t\in [0,T] \\
p(T)=-G_T,
         \end{array}
\right.
\end{equation}
then 
\begin{equation}\label{ec5.2}
\int_{\mathbb{R}^d}\int_0^T{\cal K}(\xi )\cdot b^*(\rho ^{\zeta ^*})\nabla p \ dt \ dx
-\int_{\mathbb{R}^d}\xi (x){{\partial Q}\over {\partial z}}(x,\zeta ^*(x))dx\leq 0,
\end{equation}
for any $\xi \in L^{\infty }(\mathbb{R}^d)$ such that
$\zeta ^*+\varepsilon \xi\in {\cal M}$, for any $\varepsilon >0$ sufficiently small.
\vspace{2mm}

\noindent
Proof. \rm We shall use the results obtained for problem $(P_h)$ in the previous section.

For the sake of clarity we shall renote by $\zeta _h$ (instead of $\zeta ^*_h$) an optimal control for $(P_h)$, and assume that $G$ is time-independent.
Let 
$$p_h^{\zeta _h}(t,x)=\left\{ \begin{array}{ll}
-G_T(x), \quad &\mbox{\rm if } t=T \\
p_h^{\zeta _h,i-1}(x), \quad &\mbox{\rm if } t\in [(i-1)h,ih) , \ i=1,2,...,N.
\end{array}
\right. $$

Multiplying (\ref{eq4.1}) by $-\Delta p_h^{\zeta _h,i}$ we get that
$$\int_{\mathbb{R}^d}|\nabla p_h^{\zeta _h,i}|^2_ddx +h\int_{\mathbb{R}^d}\beta '(\rho _h^{\zeta _h,i+1})|\Delta p_h^{\zeta _h,i}|^2dx =\int_{\mathbb{R}^d}\nabla p_h^{\zeta _h,i}\cdot \nabla p_h^{\zeta _h,i+1}dx$$
$$+h\int_{\mathbb{R}^d}G(x,\zeta _h(x))\Delta p_h^{\zeta ,i}dx-h\int_{\mathbb{R}^d}{\cal K}(\zeta _h)(b^*)'(\rho _h^{\zeta _h,i+1})\cdot \nabla p_h^{\zeta _h,i}\Delta  p_h^{\zeta _h,i}dx,$$
and so there exist two positive constants $M_1^0, \ M_2^0$ such that
$${1\over 2}\int_{\mathbb{R}^d}|\nabla p_h^{\zeta _h,i}|^2_ddx +{{h\gamma _0}\over 2}\int_{\mathbb{R}^d}|\Delta p_h^{\zeta _h,i}|^2dx \leq {1\over 2}\int_{\mathbb{R}^d}|\nabla p_h^{\zeta _h,i+1}|^2_ddx$$
$$+{{hM_1^0}\over 2} \int_{\mathbb{R}^d}|\nabla p_h^{\zeta _h,i}|^2_ddx+{{hM_2^0}\over 2} \int_{\mathbb{R}^d}G(x)^2dx.$$
It follows that
$$(1-hM_1^0)\int_{\mathbb{R}^d}|\nabla p_h^{\zeta _h,i}|^2_ddx +h\gamma _0\int_{\mathbb{R}^d}|\Delta p_h^{\zeta _h,i}|^2dx \leq \int_{\mathbb{R}^d}|\nabla p_h^{\zeta _h,i+1}|^2_ddx$$
$$+hM_2^0 \int_{\mathbb{R}^d}G(x)^2dx.$$
After an easy evaluation we have that for $h>0$ sufficiently small
$$(1-hM_1^0)^N\int_{\mathbb{R}^d}|\nabla p_h^{\zeta _h,i}|^2_ddx +h\gamma _0(1-hM_1^0)^N\sum_{i=0}^{N-1}\int_{\mathbb{R}^d}|\Delta p_h^{\zeta _h,i}|^2dx \leq \int_{\mathbb{R}^d}|\nabla G_T|^2_ddx$$
$$+hM_2^0 [1+(1-hM_1^0)+(1-hM_1^0)^2+\cdots ]\int_{\mathbb{R}^d}G(x)^2dx.$$
This implies that there exists a nonnegative constant $\tilde{M}$ such that for any $t\in [0,T]$:
$$\int_{\mathbb{R}^d}|\nabla p_h^{\zeta _h}(t,x)|^2_ddx +\int_0^T\int_{\mathbb{R}^d}|\Delta p_h^{\zeta _h}(t,x)|^2dx \ dt \leq \tilde{M}.$$
So, $\{ \nabla p_h^{\zeta _h}\} $ is bounded in $L^{\infty }(0,T; L^2(\mathbb{R}^d)^d)$ and $\{ \Delta p_h^{\zeta _h}\} $ is bounded in $L^2(0,T; L^2(\mathbb{R}^d))$
(for $h>0$ sufficiently small).

If we multiply now (\ref{eq4.1}) by $p_h^{\zeta _h,i}$ and integrate over $\mathbb{R}^d$ we obtain in the same manner that
$\{ p_h^{\zeta _h}\} $ is bounded in $L^{\infty }(0,T; L^2(\mathbb{R}^d))$.

We may infer that there exists a subsequence (also denoted by $\{ p_h^{\zeta _h}\} $) and a function $p$ such that
$$\left\{ \begin{array}{ll}
p_h^{\zeta _h}\longrightarrow p \quad &\mbox{\rm weakly in } L^2(0,T;L^2(\mathbb{R}^d)), \\
\nabla p_h^{\zeta _h}\longrightarrow \nabla p \quad &\mbox{\rm weakly in } L^2(0,T;L^2(\mathbb{R}^d)^d), \\
\Delta p_h^{\zeta _h}\longrightarrow \Delta p \quad &\mbox{\rm weakly in } L^2(0,T;L^2(\mathbb{R}^d)).
\end{array} 
\right. $$
By  (\ref{eq4.1}) we have that 
$$p_h^{\zeta _h,i}(x)=-G_T(x)+\int_{ih}^T\beta '(\rho _h^{\zeta _h}(s,x))\Delta p_h^{\zeta _h}(s,x)ds+\int_{ih}^T{\cal K}(\zeta _h)(x)(b^*)'(\rho _h^{\zeta _h}(s,x))\cdot \nabla p_h^{\zeta _h}(s,x)ds $$
$$-(T-ih)G(x) \quad \mbox{\rm in } L^2(\mathbb{R}^d).$$
Since $\zeta _h\longrightarrow \zeta ^*$ a.e. $x\in \mathbb{R}^d$ (on a subsequence), we may conclude that for any $t\in [0,T]$:
$$p(t)=-G_T+\int_t^T\beta '(\rho ^{\zeta ^*}(s))\Delta p(s)ds+\int_t^T{\cal K}(\zeta ^*)(b^*)'(\rho ^{\zeta ^*}(s))\cdot \nabla p(s)ds $$
$$-\int_t^TG(x)ds \quad \mbox{\rm in } L^2(\mathbb{R}^d).$$
We may conclude that $p\in W^{1,2}([0,T];L^2(\mathbb{R}^d))\cap L^2(0,T;H^2(\mathbb{R}^d))$ and that $p$ satisfies (\ref{ec5.1}).
The uniqueness of the solution to (\ref{ec5.1}) follows in a standard manner.
\vspace{2mm}

By (\ref{eq4.2}) we have that 
$$\int_{\mathbb{R}^d}\int_0^T{\cal K}(\xi )\cdot b^*(\rho _h^{\zeta ^*_h})\nabla p_h^{\zeta ^*_h}dt \ dx
-\int_{\mathbb{R}^d}\xi [{{\partial Q}\over {\partial z}}(x,\zeta _h^*(x))+(\zeta ^*_h(x)-\zeta ^*(x))]dx\leq 0,$$
for any $\xi \in L^{\infty }(\mathbb{R}^d)$ such that
$\zeta_h^*+\varepsilon \xi\in {\cal M}$, for any $\varepsilon >0$ sufficiently small, which immediately implies (\ref{ec5.2}).
\vspace{3mm}

\noindent
\it Proof of Lemma 5.2. \rm
\vspace{2mm}
We have that
$$\rho _h^{\zeta _h,i+1}-h\Delta \beta (\rho _h^{\zeta _h,i+1})+h\nabla \cdot ( {\cal K}(\zeta _h) b^*(\rho _h^{\zeta _h,i+1}))=\rho ^{\zeta _h,i}_h$$
and
$$\rho _h^{\zeta ,i+1}-h\Delta \beta (\rho _h^{\zeta ,i+1})+h\nabla \cdot ( {\cal K}(\zeta _h)  b^*(\rho _h^{\zeta ,i+1}))=\rho ^{\zeta ,i}_h
+h\nabla \cdot (( {\cal K}(\zeta _h) - {\cal K}(\zeta ) )b^*(\rho _h^{\zeta ,i+1})).$$
Using the properties of $A^{\zeta _h}$, $A^{\zeta }$ we get that
$$\| \rho_h^{\zeta _h,i+1}-\rho _h^{\zeta ,i+1}\| _{L^1}\leq \| \rho ^{\zeta _h,i}_h-\rho ^{\zeta ,i}_h
-h\nabla \cdot (( {\cal K}(\zeta _h) - {\cal K}(\zeta ) )b^*(\rho _h^{\zeta ,i+1}))\|_{L^1}$$
$$\leq  \| \rho ^{\zeta _h,i}_h-\rho ^{\zeta ,i}_h\| _{L^1}+h\|\nabla \cdot ( {\cal K}(\zeta _h) - {\cal K}(\zeta ) )b^*(\rho _h^{\zeta ,i+1})\|_{L^1}
+h\| ( {\cal K}(\zeta _h) - {\cal K}(\zeta )(b^*)'(\rho _h^{\zeta ,i+1})\cdot \nabla \rho_h^{\zeta ,i+1} \|_{L^1}.$$
By summation we obtain that for any $t\in [0,T]$:
$$\| \rho_h^{\zeta _h}(t)-\rho _h^{\zeta }(t)\| _{L^1}\leq \| |\nabla \cdot ( {\cal K}(\zeta _h) - {\cal K}(\zeta ) )|\int_0^T|b^*(\rho _h^{\zeta }(s))|ds\| _{L^1}$$
$$+\| | {\cal K}(\zeta _h) - {\cal K}(\zeta ) |_d\int_0^T|(b^*)'(\rho _h^{\zeta }(s))\cdot \nabla \rho_h^{\zeta }(s)|_dds\|_{L^1}$$
$$\leq \tilde{M}_1\| {\cal K}(\zeta _h) - {\cal K}(\zeta ) \|_{L^2}+\tilde{M}_2\|\nabla \cdot ( {\cal K}(\zeta _h) - {\cal K}(\zeta ) )\|_{L^2},$$
where $\tilde{M}_1, \ \tilde{M}_2$ are positive constants. Using now Lemma 5.1 we may infer that
$$\rho_h^{\zeta _h}(t)-\rho_h^{\zeta }(t)\longrightarrow 0 \quad \mbox{\rm in } L^1(\mathbb{R}^d), \ \mbox{\rm uniformly for } t\in [0,T].$$
On the other hand
$$\rho_h^{\zeta }(t)-\rho^{\zeta }(t)\longrightarrow 0 \quad \mbox{\rm in } L^1(\mathbb{R}^d), \ \mbox{\rm uniformly for } t\in [0,T]$$
via Crandall-Liggett theorem. We may infer that
$$\rho_h^{\zeta _h}(t)\longrightarrow \rho^{\zeta }(t) \quad \mbox{\rm in } L^1(\mathbb{R}^d), \ \mbox{\rm uniformly for } t\in [0,T].$$

\vspace{5mm}

\noindent 
\bf Remark 5.1. \rm By (\ref{ec5.2}) we get that for almost any $x\in B(0_d;R_0)$ we have that
$$\zeta ^*(x)=\left\{ \begin{array}{ll}0, & \mbox{\rm if } \displaystyle \int_{\mathbb{R}^d}\nabla {\cal G}_R(x'-x)\cdot \int_0^T[b^*(\rho ^{\zeta ^*})\nabla p](t,x')dt \ dx'+{{\partial Q}\over {\partial z}}(x,\zeta ^*(x))>0 \\
~  & ~ \\
\tilde{M}_0,  & \mbox{\rm if } \displaystyle \int_{\mathbb{R}^d}\nabla {\cal G}_R(x'-x)\cdot \int_0^T[b^*(\rho ^{\zeta ^*})\nabla p](t,x')dt \ dx'+{{\partial Q}\over {\partial z}}(x,\zeta ^*(x))<0  .
\end{array}
\right. $$

\section{Final comments}

Problem $(P)$ has an obvious meaning if we consider $G$ and $G_T$ to be approximations of the characteristic functions of two compact subsets $K\subset [0,T]\times \mathbb{R}^d$ and $K_T\subset \mathbb{R}^d$.
\vspace{3mm}

\noindent
\bf Remark 6.1. \rm The results in the present paper can be easily extended to the case of 
$$\zeta \in {\cal M}_0=\{ \theta \in L^{\infty }(\mathbb{R}^d); \ 0\leq \theta (x)\leq \tilde{M}_0 \ \mbox{\rm a.e. } x\in \mathbb{R}^d, \ \theta (x)=0 \ \mbox{\rm a.e. } x\in Ext (D)\} ,$$
where $D$ is a bounded and open subset of $\mathbb{R}^d$.
\vspace{3mm}

\noindent
\bf Remark 6.2. \rm Due to Lemma 5.1 it follows that the results in the previous sections 3-5 hold if we consider for $G$ the weaker hypotheses
$$G\in L^2((0,T)\times \mathbb{R}^d), \quad 0\leq G(t,x) \ \mbox{\rm a.e. } (t,x)\in (0,T)\times \mathbb{R}^d.$$
\vspace{3mm}

Let us briefly discuss the situation when $Q\equiv 0$ (i.e. the cost of the control $\zeta$ is negligible). Problem $(P)$ becomes
$$\underset{\zeta \in {\cal M}}{\text{Minimize}} \  \int_0^T \int_{\mathbb{R}^d}G(t,x)\rho ^{\zeta}(t,x) dx \ dt+\int_{\mathbb{R}^d}G_T(x)\rho ^{\zeta}(T,x)dx,\leqno({\bf \tilde{P}})$$
where $\rho ^{\zeta}$ is the weak solution to the following nonlinear Fokker-Planck equation (\ref{nfp1}).
For any $h={T\over N}$  (where $N\in \mathbb{N}^*$) sufficiently small we consider the following approximating optimal control problem 
$$\underset{\zeta \in {\cal M}}{\text{Minimize}} \  \int_0^T \int_{\mathbb{R}^d}G(t,x)\rho ^{\zeta }_h(t,x) dx \ dt+\int_{\mathbb{R}^d}G_T(x)\rho ^{\zeta }_h(T,x)dx, \leqno({\bf \tilde{P}_h})$$
where $\rho ^{\zeta }_h$ and $\rho ^{\zeta ,i}_h$ were defined in section 2.

Since problem $(\tilde{P})$ is a particular case of problem $(P)$, it follows that there exists at least one optimal control for $(\tilde{P})$. The existence of an optimal control $\tilde{\zeta }_h^*$ for problem $(\tilde{P_h})$ is proven in the same way as for problem $(P_h)$.
\vspace{2mm}

The following results follow in an analogous manner to Theorem 5.3 and Theorem 4.5.
\vspace{3mm}

\noindent
\bf Theorem 6.1 (The relationship between $\bf (\tilde{P}_h)$ and $\bf (\tilde{P})$). \it If $\tilde{\zeta }^*$ is a weak * accumulation point in $L^{\infty}(\mathbb{R}^d)$ for $\{ \tilde{\zeta }_h^*\}$
(and denote the convergent subsequence by $\{ \tilde{\zeta }_h^*\}$ as well), then
\begin{itemize}
\item[(l)] $\tilde{\zeta }^*$ is an optimal control for problem $(\tilde{P})$;
\item[(ll)] $\tilde{I}_h(\tilde{\zeta }_h^*)\longrightarrow \tilde{I}(\tilde{\zeta }^*)$;
\item[(lll)] $\tilde{I}(\tilde{\zeta }_h^*)\longrightarrow \tilde{I}(\tilde{\zeta }^*)$, as $h\rightarrow 0$.
\end{itemize}
Here $\tilde{I}$ and $\tilde{I}_h$ are the cost functionals for problems $(\tilde{P})$ and $(\tilde{P}_h)$, respectively.
\vspace{3mm}
\rm

If we replace the assumption $G_T\in C_b(\mathbb{R}^d)\cap H^1(\mathbb{R}^d)$ by the weaker one 
$$G_T\in C_b(\mathbb{R}^d)\cap L^2(\mathbb{R}^d),$$
then the next result holds

\noindent
\bf Theorem 6.2. \it If $p_h^{\tilde{\zeta }^*_h,i+1}$ is the unique weak solution to $(\ref{eq4.1})$
corresponding to $\zeta :=\tilde{\zeta }^*_h$, then
\begin{equation*}
\int_{\mathbb{R}^d}{\cal K}(\xi )\cdot [\sum_{i=1}^N[hb^*(\rho _h^{\tilde{\zeta }^*_h,i})\nabla p_h^{\tilde{\zeta }^*_h,i-1}]dx\leq 0,
\end{equation*}
for any $\xi \in L^{\infty }(\mathbb{R}^d)$ such that
$\tilde{\zeta }_h^*+\varepsilon \xi\in {\cal M}$, for any $\varepsilon >0$ sufficiently small.

\rm

\vspace{3mm}

\noindent
\bf Remark 6.3. \rm We can adapt the proofs in the previous sections to investigate the more general optimal control problem
$$\underset{\zeta \in {\cal M}}{\text{Minimize}} \  \int_0^T \int_{\mathbb{R}^d}G(t,x,\rho ^{\zeta}(t,x)) dx \ dt+\int_{\mathbb{R}^d}G_T(x,\rho ^{\zeta}(T,x))dx+\int_{\mathbb{R}^d}Q(x,\zeta (x))dx,\leqno({\bf P_0})$$
where $\rho ^{\zeta}$ is the weak solution to (\ref{nfp1}) if appropriate assumptions on $G, \ G_T$ are considered.

\end{document}